\documentclass{tac}
\usepackage{xy}
\usepackage{MnSymbol}
\usepackage{tikz}

\xyoption{all}

\def\mathcaldef#1{\expandafter\def\csname#1\endcsname{{\cal#1}}}

\def\"{``}
\def\q{\quad}
\def\qq{\quad\quad}
\def\qv{\qq ;\qq}

\def\iso{\,\cong\,}
\def\la{\langle}
\def\ra{\rangle}
\def\adj{\dashv}
\def\op{^{\rm op}}
\def\ov{\overline}
\def\ul{\underline}
\def\tm{\times}
\def\otm{\otimes}
\def\lam{\lambda}
\def\t{_\blacktriangleright}

\def\adots{a_1,\cdots , a_n}
\def\Xdots{X_1,\cdots , X_n}
\def\ldots{\lam_1,\cdots ,\lam_n}
\def\ladots{\lam_1 a_1+\cdots +\lam_n a_n}

\mathcaldef{C}
\mathcaldef{D}
\mathcaldef{M}
\mathcaldef{N}
\mathcaldef{S}
\mathcaldef{P}
\mathcaldef{V}

\mathbfdef{Set}
\mathbfdef{Pos}
\mathbfdef{cMon}
\mathbfdef{cComon}
\mathbfdef{Mon}
\mathbfdef{Mod}
\mathbfdef{Mlt}
\mathbfdef{fpMlt}
\mathbfdef{sMlt}
\mathbfdef{Rep}
\mathbfdef{fpRep}
\mathbfdef{Cat}
\mathbfdef{fpCat}
\mathbfdef{Sum}
\mathbfdef{fpSum}
\mathbfdef{mltfb}
\mathbfdef{mltdf}

\mathrmdef{id}
\mathrmdef{I}
\mathrmdef{m}
\mathrmdef{v}
\mathrmdef{s}
\mathrmdef{d}

\def\obj{{\rm obj\,}}

\newtheorem{prop}{Proposition}
\newtheorem{corol}{Corollary}
\let\pf\proof
\let\epf\endproof
\def\eq{\begin{equation}}
\def\eeq{\end{equation}}

\author{Claudio Pisani}
\address{via Saluzzo 67,\\ 10125 Torino, Italy.}
\title{Some remarks on multicategories \\ and additive categories}
\copyrightyear{2013}
\eaddress{pisclau@yahoo.it}

\begin{document}

\maketitle

\begin{abstract}

Categories are coreflectively embedded in multicategories via the \"discrete cocone" 
construction, the right adjoint being given by the monoid construction.
Furthermore, the adjunction lifts to the \"cartesian level": preadditive categories  
are coreflectively embedded (as theories for many-sorted modules)
in cartesian multicategories (general algebraic theories).
In particular, one gets a direct link between two ways of considering modules over a rig,
namely as additive functors valued in commutative monoids or as models of the theory 
generated by the rig itself.


\end{abstract}

\section{Introduction}
\label{intro}

The present work develops around the idea that a (semi)module over a rig $\C$
is a model of the theory whose operations are finite sequences of arrows in $\C$,
giving linear combinations:
\eq   \label{intro5}
\la \ldots \ra : \la \adots \ra \mapsto \ladots \qv  \la\ra : \la\ra \mapsto 0 
\eeq
(Recall that a rig, a ring without negatives, is a monoid enriched in commutative monoids.)
Indeed, the rig $\C$ generates a one-object cartesian multicategory $\C\t$
such that modules over $\C$ are morphisms of cartesian multicategories $\C\t \to \Set$,
as we are going to illustrate below in this introduction.

The relation between monoidal categories and multicategories is well-known 
(see for instance \cite{lambek}, \cite{hermida} and \cite{leinster});  
in particular there is a functor $(-)_\otm : \Cat^\otm \to \Mlt$
from the category of monoidal categories and lax monoidal functors
to the category of multicategories and (multi)functors, with
\eq  \label{}
\M_\otm(\Xdots;X) \iso \M(X_1\otm \cdots \otm X_n;X)
\eeq
that restricts to an equivalence $\Cat^\otm \simeq \Rep$,
where $\Rep \subset \Mlt$ is the full subcategory of representable multicategories.

Less known is the important fact that also categories (with functors) are embedded
in $\Mlt$ via the discrete cocone functor $(-)\t:\Cat\to\Mlt $:
\[ 
\C\t(\Xdots;X) \iso \C(X_1;X)\tm \cdots \tm \C(X_n;X)
\]
So $n$-arrows in $\C\t$ are sequences $\la \ldots \ra$ of concurrent arrows in $\C$
with the obvious composition, e.g. 
 \[
\la\lam,\mu,\nu\ra(\la \alpha,\beta\ra,\la\,\ra,\la \gamma \ra) = \la \lam\alpha,\lam\beta,\nu\gamma \ra 
\]
and for any functor $F:\C\to\D$
\[
F\t \la \ldots \ra = \la F\lam_1,\cdots, F\lam_n \ra
\]

This functor (whose existence I learned from \cite{hermida}; see also \cite{leinster}) 
is pivotal in our work, although it does not seem to have received much attention in the literature. 
We will denote by $\Cat\t\subset\Mlt$ the full subcategory of those
multicategories isomorphic to $\C\t$ for some category $\C$,
so that there is an equivalence $\Cat\t \simeq \Cat$.

Note that $\Sum := \Cat\t \cap \Rep \subset \Mlt$ is equivalent to the category
of categories with finite coproducts (and all functors);
indeed, $\C\t$ is representable iff $\C$ has finite coproducts.
Note also that multicategories of the form $\M_\otm$ and $\C\t$ give two main instances
of promonoidal categories \cite{day}.

The notion of cartesian multicategory we refer to is that in \cite{fiore2}, which actually
does not even mention explicitly multicategories (see however \cite{gould}).
It is essentially a multicategory $\M$ with an action of mappings (over $\obj\M$) 
${\bf n} = \{1,\cdots,n\}\to {\bf m} = \{1,\cdots,m\}$  
on its hom-sets which is compatible with composition (see Section \ref{cart}).
We denote by $\fpMlt$ the category of cartesian multicategories and their 
morphisms, with the obvious underlying functor to (plain) multicategories $U:\fpMlt \to \Mlt$.
We will sometimes refer to cartesian multicategories and their morphisms as 
multicategories with an \"fp-structure" and \"fp-functors" respectively.

So as monoidal categories and categories give two main instances of multicategories, 
finite product categories and preadditive categories (that is enriched in commutative monoids) 
give rise, for somewhat opposite reasons that we will illustrate presently, 
to two main instances of cartesian multicategories. 
More is true: the full subcategory $\fpRep\subset\fpMlt$ above (with respect to $U$) $\Rep$ 
is equivalent to the category of finite product categories and finite product preserving functors
while the full subcategory $\fpCat\t\subset\fpMlt$ above $\Cat\t$ 
is equivalent to the category of preadditive categories and additive functors
(propositions \ref{cart4} and \ref{cart5}).
Then $\fpSum = \fpCat\t \cap \fpRep \subset \fpMlt$ 
is equivalent to the category of additive categories; see the figure on page 14.

Let us come back to our one-object case.
If $\C$ is a rig, the operad $\C\t$ has the following fp-structure: 
bijections $\bf n\to n$ act obviously on sequences $\la \ldots \ra$, 
while the effect of monotone mappings is to sum pairs or sequences in $\C$
and to insert zeros.
For instance, if 
\eq   \label{intro7}
p: \{ 1, 2, 3 \} \to \{ 1, 2, 3 \} \qv 1\mapsto 1,\q 2\mapsto 1,\q 3\mapsto 3
\eeq
then 
\[
p\la\alpha,\beta,\gamma\ra = \la\alpha+\beta,0,\gamma\ra
\]

On the other hand, the finite product category $\Set$ gives rise to the 
cartesian multicategory $\Set_\tm$ with
\[
\Set_\tm(\Xdots ; X) = \Set(X_1\tm \cdots \tm X_n ; X)
\] 
and with actions of mappings $\bf n \to m$ given by composition with diagonals and projections;
for instance for $p: \{ 1, 2, 3 \} \to \{ 1, 2, 3 \}$ as in (\ref{intro7}) 
and for $f:X\tm X\tm X \to Y$ in $\Set$ we have
\[
p f(x,y,z) = f(x,x,z)
\]

Now, to give a morphism $\ov{( -)}:\C\t \to \Set_\tm$ in $\fpMlt$ amounts to 
\begin{enumerate}
\item
giving a set $X$ with a commutative monoid structure;
indeed, if we define 
\[
\m_n = \ov{\la 1, \cdots ,1\ra}
\] 
with $n$ 1's, then $\m_2:X\tm X \to X$ and $\m_0:1\to X$ clearly define a monoid
(since $\ov{( -)}$ is a morphism of multicategories) which is commutative 
(since it also preserves the fp-structure relative to bijections);
\item
giving (functorially) a monoid endomorphism $\ov\lam :X\to X$ for any $\lam \in\C$;
this follows from the equalities in $\C\t$
\[
\la 1,1\ra (\la\lam\ra,\la\lam\ra) = \la \lam \ra (\la 1, 1 \ra ) \qv \la\,\ra = \la \lam \ra ( \la\,\ra )
\]
which become in $\Set$ the monoid morphism conditions 
\[
\m_2 (\ov\lam,\ov\lam) = \ov\lam \m_2 \qv \m_0 = \ov\lam \m_0
\]
Note that from the equality $\la \alpha, \beta \ra = \la 1,1\ra (\la\alpha\ra,\la\beta\ra)$
in $\C\t$ one gets, using the sum notation for the monoid structure on $X$,
\[
\ov{\la \alpha, \beta\ra}(x,y) = \ov\alpha x + \ov\beta y    \qv  \ov{\la \, \ra} = 0
\]
that is the arrows in $\C\t$ become indeed \"linear combinations" in $X$;
\item
such that the fp-structure is preserved also for the non-bijective mappings; 
for instance, if $p$ is as in (\ref{intro7})
\[
p \ov{\la\alpha,\beta,\gamma\ra} = \ov{p\la\alpha,\beta,\gamma\ra} = \ov{\la\alpha +\beta,0,\gamma\ra}
\]
that is
\[
\ov\alpha x +\ov\beta x+\ov\gamma z  = \ov{(\alpha +\beta)}x+\ov 0 y +\ov\gamma z
\]
\end{enumerate} 
Summarizing, to give a morphism $\ov{( -)}:\C\t \to \Set_\tm$ in $\fpMlt$ amounts to
give a commutative monoid $(X,+,0)$ in $\Set$ which has a module structure on $\C$:
\[ \ov 1 x = x \qv \ov{\alpha\beta} x = \ov\alpha (\ov\beta x) \]
\[ \ov \alpha 0 = 0 \qv \ov\alpha (x+y) = \ov\alpha x + \ov\alpha y \]
\[ \ov 0 x = 0 \qv  \ov{(\alpha +\beta)}x = \ov\alpha x +\ov\beta x \]

Thus, for a rig $\C$, we have a correspondence
\eq     \label{intro2}
\begin{array}{c}
\fpMlt(\C\t, \Set) \\ \hline
\C-\Mod \\ \hline
\cMon\Cat(\C,\cMon(\Set)) 
\end{array}
\eeq
which is in fact just an instance of an adjunction (see Proposition \ref{cart2}):
\eq     \label{intro1}
(-)\t \adj \cMon(-) : \fpMlt \to \cMon\Cat 
\eeq


\subsection{Summary}

In Section \ref{plain} we present the tensored structure of the 2-category $\Mlt$ of multicategories,
(multi)functors and natural transformations, which implies the \"plain version" of the adjunction (\ref{intro1}): 
\eq     \label{intro3}
(-)\t \adj \Mon(-) : \Mlt \to \Cat 
\eeq

In Section \ref{cart} we move to the \"cartesian level", discussing the relationships between
cartesian multicategories, (pre)additive categories and finite product categories.

In Section \ref{fib} we briefly consider the results from a fibrational point of view:
from (\ref{intro3}) we get in particular that (strong) indexed monoidal categories
correspond to functors $\C\t \to \Cat_\tm$ in $\Mlt$, and thus (ignoring pseudo-issues for simplicity) 
also to multicategory fibrations (in the sense of \cite{hermida2}) over $\C\t$:
\eq     \label{intro4}
\begin{array}{c}
\Cat(\C,\Mon(\Cat_\tm)) \\ \hline
\Mlt(\C\t, \Cat_\tm)  \\ \hline
\mltfb/\C\t
\end{array}
\eeq

We conclude by presenting some elementary examples of multicategories arising in this way 
as the domain of a fibration over $\C\t$ for an indexed monoidal category $\C \to \Mon(\Cat)$.


The aim of the present paper is mainly expository, 
so that some proofs are omitted or just sketched;  
for a specific but not central fact (Proposition \ref{cart5}) I have to refer to \"folklore".
A more detailed account will follow.

\subsection{Terminology and notations}

Following other authors, we omit the prefix \"semi" which is sometimes 
used to refer to commutative monoids (rather than abelian groups)
enrichments.
Thus, \"preadditive" categories and \"additive" functors are $\cMon$-enriched 
categories and functors, while \"additive" categories are those preadditive
categories with biproducts (again with additive functors, that necessarily preserve
biproducts).

As already seen, we use the following notation which may be a little confusing:
$\fpCat\t$ is the category of those cartesian multicategories whose underlying 
multicategory has the form $\C\t$ (with all fp-functors).
It is equivalent to the category of preadditive categories.
On the other hand, $\fpRep$ is the category of representable cartesian multicategories
and is equivalent to $\Cat^\tm$, the category of finite product categories
and finite product preserving functors. 
In fact, it is this equivalence that suggests to extend the \"fp" prefix (for \"finite product")
to all the \"cartesian level". A similar notation is used in \cite{gould}.

So the prefix \"fp" (and the term \"cartesian" itself) refer to multicategories rather than to categories,
allowing us to emphasize the correspondence between the cartesian and the plain levels
(see the figure on page 14).


\section{The plain level}
\label{plain}

We assume that the reader is familiar with the 2-category $\Mlt$ of multicategories,
functors and natural transformations; \cite{leinster} is a good reference.


The discrete cocone functor $(-)\t:\Cat\to\Mlt$
mentioned in the introduction is in fact a full and faithful 2-functor.
For fullness, let $F:\C\t \to \D\t$ be a functor in $\Mlt$,
$\la\lam,\mu\ra \in \C\t(X,Y;Z)$, and $F \la\lam,\mu\ra = \la\lam',\mu'\ra \in \D\t(FX,FY;FZ)$;
then
\[
\lam' = F \la\lam,\mu\ra (\id_{FX},\la\,\ra) = F \la\lam,\mu\ra (F\id_X,F\la\,\ra) = 
F (\la\lam,\mu\ra (\id_X,\la\,\ra)) = F\lam
\]
so that in fact $F = (\ul F)\t$, where $\ul F : \C \to \D$ is the \"underlying" functor.

It is also easy to see that $(-)\t$ preserves products. In particular, $1\t$ is terminal in $\Mlt$.

The next proposition says that the 2-category $\Mlt$ is tensored and that
the tensor of $\M\in\Mlt$ by $\C\in\Cat$ is provided by the cartesian product
in $\Mlt$ of $\M$ and $\C\t$:
\begin{prop}
 There are natural isomorphisms
\eq     \label{}
\begin{array}{c}
\Cat(\C,\Mlt(\M,\N)) \\ \hline
\Mlt(\C\t\tm\M, \N) 
\end{array}
\eeq
\end{prop}
\pf
Straightforward calculations. 
\epf

Thus, for any $\M\in\Mlt$ we have an adjunction
\eq     \label{plain2}
(-)\t \tm \M \adj \Mlt(\M,-):\Mlt\to\Cat
\eeq
The most important particular case is obtained when $\M = 1\t$;
recall that, for $\N\in\Mlt$, the category $\Mlt(1\t,\N)$ can be identified
as the category $\Mon(\N)$ of monoids in $\N$:
\begin{corol}  \label{plain1}
The monoid functor is right adjoint to the discrete cocone functor:
\eq     \label{plain4}
(-)\t \adj \Mon(-):\Mlt\to\Cat
\eeq
\epf
\end{corol} 

\begin{remark}
Then we could say that $\C\t$ is the universal multicategory containing 
a $\C$-indexed monoid.
\end{remark}

\begin{remark}
The counit of the adjunction takes an arrow $\la \ldots \ra : \la \Xdots \ra \to X$ in $\Mon(\M)\t$
to $\m_n^X(\ldots)$ in $\M$.
\end{remark}

\begin{remark}  \label{plain3}
Since $(-)\t$ is fully faithful, $\Mon(\C\t) \iso \C$.
When $\C$ has finite sums we find again the classical fact that in the monoidal category
$(\C,+,0)$ every object carries a unique monoid structure.
The converse also holds; see for instance \cite{heunen}.
\end{remark}

\begin{corol}  
$\Cat\t$ is coreflective in $\Mlt$.
\epf
\end{corol}

\begin{remark}
This is related to an old result in \cite{fox}:
the commutative monoid construction provides a universal way to endowe a symmetric 
monoidal category with finite sums.
Indeed, an intermediate level between cartesian and plain multicategories is that of 
symmetric multicategories (see \cite{leinster}, \cite{lambek} and \cite{gould}):
multicategories in $\sMlt$ are defined as cartesian multicategories except that only 
bijective mappings are supposed to act on the hom-sets.
Any $\C\t$ has an ovious s-structure and it is easy to see that there is a symmetric
version of (\ref{plain4}):
\eq    
(-)\t \adj \cMon(-):\sMlt\to\Cat
\eeq
so that $\Cat\t$ is coreflective in $\sMlt$.
Thus, Fox's result can be seen as stating that the latter restricts to representable symmetric 
multicategories: finite sums categories are coreflective in symmetric monoidal categories.
\end{remark}

Another particular case of (\ref{plain2}) is obtained when $\M = 1$, the multicategory with just an arrow; 
it is immediate to observe that $\Mlt(1,\N)$ is the \"underlying" category $\ul\N$ of $\N$,
with $\ul\N(X;Y) = \N(X;Y)$, while $\C\t \tm 1$ is the \"linear" multicategory $\C_!$,
with $\C_!(X;Y) = \C(X;Y)$ and $\C_!(\Xdots;Y) = \emptyset$ for $n\neq 1$:

\begin{corol}  
The underlying functor is right adjoint to the linear functor:
\[ 
(-)_! \adj \ul{(-)}:\Mlt\to\Cat
\]
\epf
\end{corol}

Since $(-)_! :\Cat\to\Mlt$ is also clearly a full and faithful 2-functor,
$\Cat$ is coreflectively embedded in two ways in $\Mlt$.
Anyway, in the present work the prominent role is played by the discrete cocone embedding.

Note that $\ul{(-)}\circ (-)\t : \Cat \to \Cat$ is the identity, so that both $\ul{(-)}$ and $\Mon(-)$
provide left inverses for the discrete cocone functor.


\section{The cartesian level}
\label{cart}

\subsection{Cartesian multicategories}

The notion of cartesian multicategory has appeared in various guises (and names) in the literature;
see e.g. \cite{lambek}, where they are called \"Gentzen multicategories".
The idea is that they are to finite product categories what multicategories are to monoidal categories.
So they should fulfill at least the condition that representable cartesian multicategories are 
finite product categories, that is the tensor product is cartesian.
An elegant definition is given in \cite{fiore2} (without mentioning explicitly multicategories)
and taken over in \cite{gould} for operads.
We report here the essence of that notion.

Let $\bf N$ the full subcategory of $\Set$ which has objects ${\bf 0, 1, 2}, \cdots $, with
${\bf n} = \{ 1, 2, \cdots , n \}$ (so that ${\bf 0} = \emptyset$) and consider, 
for a multicategory $\M$, the obvious comma category ${\bf N}/\obj\M$.
 
Now, the domain of an $n$-arrow in $\M$ is in fact a mapping $\alpha:{\bf n} \to \obj\M$,
that is an object of ${\bf N}/\obj\M$; thus, for any fixed codomain object $X\in\M$
we have a mapping
\[
\obj ({\bf N}/\obj\M) \to \obj\Set  \qv  \alpha \mapsto \M(\alpha,X) = \M(\Xdots, X)
\]
To give an fp-structure on $\M$ means to extend these mappings to functors
\[
(-)_X : {\bf N}/\obj\M \to \Set
\] 
in a way compatible with composition.

For instance, let $\obj\M = \{ A, B, C, D, \cdots \}$ and let $\alpha: {\bf 3}\to\obj\M$ 
and $\beta: {\bf 3}\to\obj\M$ be the objects of ${\bf N}/\obj\M$ given by 
\[ 
\alpha: \qq 1\mapsto A,\q 2\mapsto B,\q 3\mapsto A   \qv 
\beta: \qq   1\mapsto B, \q 2\mapsto A,\q 3\mapsto C
\]
If $p:\alpha \to \beta$ in ${\bf N}/\obj\M$ is given by 
\[
p: \qq 1\mapsto 2,\q 2\mapsto 1,\q 3\mapsto 2
\]
then an fp-structure on $\M$ gives in particular a mapping
\[
p_D : \M(A, B, A;D) \to \M(B, A, C ; D)
\]
(We will hereafter omit the subscript indicating the codomain object.)

As for the compatibility conditions with respect to composition, there are two of them.
The first one is pretty obvious: 
\[ 
f(p_1 f_1,\cdots, p_n f_n) = (p_1+\cdots+ p_n)f(f_1,\cdots, f_n)
\]
that is to compose $f$ with arrows $f_i$ acted upon by $p_i$ is the same as 
composing $f$ with the $f_i$ and then acting on it with the obvious \"sum" 
of maps in ${\bf N}/\obj\M$.

The second compatibility condition concerns the case when it is $f$ 
that is acted upon by a mapping $p$:
\[ 
(pf)(f_1,\cdots, f_n) = p'(f(f_{p1},\cdots, f_{pn}))
\]
where $p'$ is a suitably defined map in ${\bf N}/\obj\M$.

Let us illustrate it by referring to the above instance.
Let 
\[
f:\la A,B,A\ra \to D \q ;\q a:\la A_1,A_2\ra\to A \q ; \q b:\la B_1,B_2\ra\to B \q ;\q c:\la C_1,C_2\ra\to C
\]
be arrows in $\M$ and let $p$ be defined as above.
Then 
\[
(pf)(b,a,c)= p'(f(a,b,a)) : \la B_1,B_2,A_1,A_2,C_1,C_2\ra \to D 
\]
where $f(a,b,a) : \la A_1,A_2,B_1,B_2,A_1,A_2 \ra\to D$ and $p':\bf \gamma \to \delta$ 
(for the obvious $\gamma$ and $\delta$) is given by the following mapping $\bf 6 \to 6$:
\[
p': \qq 1\mapsto 3,\q 2\mapsto 4,\q 3\mapsto 1, \q 4\mapsto 2,\q 5\mapsto 3,\q 6\mapsto 4
\]

The key cases are in fact those of the mappings $p$ over $\v_n : \bf n \to 1$ (contractions or diagonals),
over $\d_{1,2}:\bf 1 \to 2$ (weakenings or projections) and $\s : \bf 2 \to 2$, 
or more generally bijections $\bf n \to n$ (exchange):
\[
\v_2:\M(X,X;Z) \to \M(X;Z)  \qv \v_0:\M(\,\, ;Z) \to \M(X;Z)  
\]
\[
\d_1:\M(X;Z) \to \M(X,Y;Z) \qv \d_2:\M(Y;Z) \to \M(X,Y;Z) 
\]
\[
\s:\M(X,Y;Z) \to \M(Y,X;Z)
\]
(Note that ${\bf 0}\to \obj\M$ is initial in ${\bf N}/\obj\M$, while $A: {\bf1}\to \obj\M$
receives exactly one map from the constantly $A$-valued object ${\bf n}\to\obj\M$
and none from the other ones.)

A functor $\M\to\M'$ induces an obvious functor  ${\bf N}/\obj\M \to {\bf N}/\obj\M'$
and maps of cartesian multicategories are of course those functors which 
commute with the actions of ${\bf N}/\obj\M$ and ${\bf N}/\obj\M'$.
We so obtain the 2-category $\fpMlt$ of cartesian multicategories with the obvious
forgetful functor $U:\fpMlt \to \Mlt$.

\begin{remark}
Cartesian multicategories can be seen as (many-sorted) algebraic theories 
with the same expressive power as (many-sorted) Lawvere theories.
Indeed, there is an equivalence between cartesian multicategories and Lawvere theories
(\cite{fiore2}; see also \cite{lambek}).
In this perspective, fp-functors $\M \to \Set$ (or more generally $\M \to \N$) are the ($\N$-valued)
\"algebras" or \"models" of $\M$. See also \cite{leinster} for a discussion of those algebraic
theories (namely the \"strongly regular" ones, whose equations are \"linear") 
which can be expressed by plain operads.
\end{remark}

We have already sketched in the introduction how cartesian multicategories arise
both from finite product categories and from preadditive categories.
Let us state this fact more formally. 
We denote by $\Cat^\tm$ the category of finite product categories and finite product preserving functors,
while $\cMon\Cat$ is the category of categories enriched in $\cMon = \cMon(\Set)$,
that is preadditive categories and additive functors.
\begin{prop}  \label{cart3}
There are fully faithful functors 
\[
(-)_\tm : \Cat^\tm \to \fpMlt   \qv  (-)\t : \cMon\Cat \to \fpMlt 
\]
\end{prop}
\pf
Given a finite product category $\C\in\Cat^\tm$ we get $\C_\tm\in\fpMlt$ with
\[ 
\C_\tm(\Xdots;X) \iso \C(X_1\tm \cdots \tm X_n;X)
\]
and with contractions and weakenings given by diagonals and projections:
\[
\v_2:\C(X\tm X;Z) \to \C(X;Z)  \qv \v_2 : f \mapsto f\circ\Delta 
\]
\[
\v_0: \C(1;Z) \to \C(X;Z) \qv  \v_0 : x \mapsto x\circ !_X
\]
\[
\d_1:\C(X;Z) \to \C(X\tm Y;Z)  \qv  \d_1 :  f \mapsto f \circ\pi_1
\]
\[
\s:\C(X\tm Y;Z) \to \C(Y\tm X;Z) \qv  \d_1 :  f \mapsto f \circ \la \pi_2,\pi_1 \ra
\]

Given a preadditive category $\C\in\cMon\Cat$ we get $\C\t\in\fpMlt$ with
\[ 
\C\t(\Xdots;X) \iso \C(X_1;X)\tm \cdots \tm \C(X_n;X)
\]
as in Section \ref{plain} and with contractions and weakenings given 
by sums and zero's insertions:
\[
\v_2:\C(X;Z)\tm \C(X;Z) \to \C(X;Z)  \qv \v_2 : \la f,g \ra  \mapsto f+g 
\]
\[
\v_0: 1 \to \C(X;Z) \qv  v_0 : * \mapsto 0_{X,Z}
\]
\[
\d_1:\C(X;Z) \to \C(X;Z) \tm \C(Y;Z)  \qv  \d_1 :   f \mapsto \la f,0 \ra
\]
\[
\s:\C(X;Z)\tm \C(Y;Z) \to \C(Y;Z)\tm\C(X;Z)  \qv \s : \la f,g \ra  \mapsto \la g,f \ra 
\]
Both the assignments are easily seen to extend to faithful functors.
For the fullness of $(-)_\tm$ note that if $F:\C_\tm \to \D_\tm$ is an fp-functor
then the underlying $\ul F : \C \to \D$ preserves projections, since
\[
\ul F \pi^{X,Y}_1 = F(\d_1 \id_X) = \d_1 \id_{FX} = \pi^{\ul FX,\ul FY}_1
\]
so that $F = (\ul F)_\tm$, with $\ul F : \C \to \D$ in $\Cat^\tm$.

As for the fullness of $(-)\t$, let $F:\C\t \to \D\t$ be an fp-functor with underlying $\ul F : \C \to \D$.
Recall from Section \ref{plain} that $F = (\ul F)\t$, that is $F\la f,g \ra = \la \ul F f, \ul F g \ra$,
so that it is enough to show that $\ul F$ is additive:
\[
\ul F(f+g) = F(\v_2 \la f,g \ra) = \v_2 F\la f,g \ra = \v_2 \la \ul F f, \ul F g \ra = \ul Ff + \ul Fg
\]
\epf


In the other direction we have:
\begin{prop}
 $\cMon(-)$ gives a functor $\fpMlt \to \cMon\Cat$.
\end{prop}
\pf
If $\M$ is a cartesian multicategory, then $\cMon(\M)$ is a preadditive category.
This can be seen as in the case $\M = \Set$ but using explicitly the fp-structure of $\M$.
In particular, if $f,g:X\to Y$ are maps in $\cMon(\M)$ their sum is given by
\[
f+g = \v_2\, m_2^Y(f,g)
\]
Again, the assignment is easily seen to extend to a functor.
\epf

\begin{prop}  \label{cart2}
There is an adjunction
\eq     \label{cart6}
(-)\t \adj \cMon(-) : \fpMlt \to \cMon\Cat 
\eeq
\end{prop}
\pf
Just follow the trace given in the introduction,
with minor modifications to take in account the many-sorted case
and the generality of $\M$ (with respect to $\Set$).
\epf

Now we characterize (up to equivalence) the codomain of the functors
of Proposition~\ref{cart3}. 
Recall that in the introduction we defined $\fpCat\t$ and $\fpRep$
as the full subcategories of $\fpMlt$ over $\Cat\t$ and $\Rep$ respectively,
so that $(-)\t$ and $(-)_\tm$ are in fact functors 
\[
(-)\t : \cMon\Cat \to \fpCat\t  \qv (-)_\tm : \Cat^\tm \to \fpRep
\]
\begin{prop}   \label{cart4}
The adjunction {\rm (\ref{cart6})} restricts to an equivalence
\[
\fpCat\t \simeq \cMon\Cat
\]
\end{prop}
\pf
If $\C\t\in\Cat\t$ has an fp-structure then contractions $\v_0$ and $\v_2$
give a monoid structure on the hom-sets of $\C$ that distributes
over compositions because of the compatibility conditions.
Thus any $\M\in\fpCat\t$ is indeed isomorphic to $\C\t$
for a preadditive $\C$.
\epf

\begin{corol}
 $\fpCat\t$ is coreflective in $\fpMlt$.
 \epf
\end{corol}

\begin{prop}   \label{cart5}
There is an equivalence
\[
\fpRep \simeq \Cat^\tm  
\]
More precisely, if a representable multicategory $\M$ has
an fp-structure then it is unique and the representing tensor
is in fact a product (so that $\M \iso \C_\tm$ for $\C\in\Cat^\tm$).
\end{prop}
\pf
As stated at the beginning of this section, the proposition should be true
for any adequate definition of cartesian multicategory.
I do not have myself a complete proof but the missing step seems to be folklore. 
First note that, for any $\M\in\fpRep$, the representing tensor product is symmetric
(because of the acting bijections, in particular the exchange $\s:{\bf 2}\to{\bf 2}$).
Furthermore, for any multicategory $\M$ the hom-sets
functors $\M^n:\ul\M\op\tm\cdots\tm\ul\M\op\tm\ul\M\to\Set$
(where $\ul\M$ is the underlying category) induce (via the diagonal
$\ul\M\op\to\ul\M\op\tm\cdots\tm\ul\M\op$) \"diagonal"
endoprofunctors $\M^n_\Delta$ on $\ul\M$:
\[
\M^n_\Delta(X,Y) = \M(X,\cdots,X;Y)
\]
Now, because of the compatibility conditions, the contraction and 
exchange mappings $\v_n$ and $\s$ are profunctor morphisms:
\[
\M^n_\Delta \to \M^1_\Delta = \ul\M  \qv \M^2_\Delta \to \M^2_\Delta
\]
with $\v_2\s = \v_2$ 
(in fact we have a functor ${\bf N} \to \Set^{\ul\M\op\tm\ul\M}$);
thus one gets (Yoneda) a natural commutative comonoid structure
\[
X\to X\otm X  \qv X\to1
\]
on each object $X\in\M$.
The missing step is then the dual of the following statement:

{\em Given a symmetric monoidal category $(\C,\otm,\I)$,  
if the forgetful functor $\cMon(\C) \to \C$ is a split epi in $\Cat$ 
then $\I$ is initial and $\otm$ is the coproduct in $\C$.}

Among the papers which cite (without proof) this fact there are 
\cite{lamarche}, \cite{fiore} and \cite{pavlovic};
an ingenious example by Jeff Egger shows that the commutativity hypothesis is necessary.
\epf

\begin{corol}
The full subcategory 
\[
\fpSum := \fpCat\t \cap \fpRep \subset \fpMlt
\] 
of those cartesian multicategories which are over $\Sum$ can be identified with both
\begin{enumerate}
\item
preadditive categories with finite sums (and additive functors);
\item
categories with finite sums and with finite products which are naturally isomorphic 
(and finite (bi)product preserving functors).
\end{enumerate}
\epf
\end{corol}

Thus, we find again the result that 
additive categories (preadditive categories with biproducts)
can be characterized both as preadditive categories with sums
and as categories with finite sums and products and with natural
isomorphisms $X + Y \to X\tm Y$ (see \cite{lack}).
Furthermore, the additive functors between additive categories
coincide with the biproduct preserving functors.

\begin{remark}
The forgetful functor $U$ has a left adjoint (see \cite{gould} for the case of operads):
\[
F \adj U : \fpMlt \to \Mlt
\]
In particular, if $\C\t\in\Cat\t$, $F\C\t$ gives the usual free preadditive category on $\C\in\Cat$.
\end{remark}


\section{Multifibrations}
\label{fib}

From the adjunction $(-)\t \adj \Mon(-):\Mlt\to\Cat$
%
of Corollary \ref{plain1}, we get in particular isomorphisms
\eq     \label{fib2}
\begin{array}{c}
\Cat(\C,\Mon(\Cat_\tm)) \\ \hline
\Mlt(\C\t, \Cat_\tm) 
\end{array}
\eeq

Now, an obvious generalization of the Grothendieck construction shows that
there is a correspondence between (covariant) fibrations of multicategories 
(see \cite{hermida2}) over $\M$ and (pseudo)functors $\M \to \Cat_\tm$: 
\eq     \label{fib3}
\begin{array}{c}
\Mlt(\M, \Cat_\tm) \\ \hline
\mltfb/\M
\end{array}
\eeq

Thus we have a correspondence
\eq     \label{fib4}
\begin{array}{c}
\Cat(\C,\Mon(\Cat_\tm)) \\ \hline
\mltfb/\C\t
\end{array}
\eeq
(more precisely, above the line we should consider pseudofunctors from $\C$
to the category representable multicategories with functors that preserve universal $n$-arrows)
which can be seen directly as follows.
 
Given $\M:\C \to \Mon(\Cat_\tm)$ and $a_i\in\M X_i$, arrows $\la\adots\ra \to a$ in 
the domain $\widehat\M$ of the corresponding fibration are sequences of arrows $\lam_i:X_i\to X$ in $\C$ 
plus a morphism $\lam_1 a_1 \otm \cdots \otm \lam_n a_n \to a$ in $\M X$.
Convesely, suppose we have a fibration over $\C\t$ and an object $X:1\to\C$; the corresponding
$X\t:1\t \to\C\t$ induces, by pullback, a fibration over $1\t$, giving a monoidal category $\M X$;
then, as in the standard case, one extends this to a (pseudo)functor $\M:\C \to \Mon(\Cat_\tm)$.

\begin{remark}
If $\C$ has finite sums, $\C\t$ is representable (that is fibered over $1\t$) so that
the domain of the fibration corresponding to an indexed monoidal category
is also representable; thus we get as a particular case of (\ref{intro4})
(the dual of) the correspondence in \cite{shulman}, for a finite product category $\C$, between
$\C\op$-indexed monoidal categories and monoidal fibrations over $\C$.
\end{remark}

Of course, (\ref{fib4}) restricts to discrete (op)fibrations:
\eq     \label{fib5}
\begin{array}{c}
\Cat(\C,\Mon(\Set_\tm)) \\ \hline
\mltdf/\C\t
\end{array}
\eeq

The idea is that a functor $\M:\C \to \Mon(\Set_\tm)$ provides a way to combine (rather than just to compare,
as in the case of functors $\C \to \Set$) objects of various sorts (elements $a_i\in\M(X_i) ,\, X_i\in\C$)
producing another object $\ladots$ (where $\lam_i\in\C$ expresses how $a_i$ is used).
The domain $\widehat\M$ of the resulting fibration is the multicategory of such combinations
and the functor $\widehat\M \to \C\t$ gives the \"abstract shape" of each concrete 
combination in $\widehat\M$.


Let us see some examples.

\subsection{Linear combinations}
\label{ex1}

A (single sorted) instance is given by the case of a module over a rig (say a vector space) 
discussed in the introduction.
In this case, arrows in $\widehat\M$ are true linear combinations and $\widehat\M(\adots;a)$
is not empty iff $a$ is spanned by the $a_i$. 
Composition in $\widehat\M$ is just the natural idea that a combination of combinations gives a combination.
The fact that $\M$ is actually a module (that is $\C$ is a rig and $\M:\C \to \cMon(\Set_\tm)$ is additive)
is reflected in the fact that $\widehat\M$ has itself an fp-structure and $\widehat\M \to \C\t$ preserves it.

\subsection{Combinations of figures}
\label{ex2}

The isometry group $\C$ of the plane $X$ acts on \"figures" $A\subseteq X$;
it in fact acts on the monoid $(\P X , \cup ,\emptyset)$.
Thus we obtain the multicategory $\widehat\M$ with figures as objects and arrows
\eq   \label{}
\la \ldots \ra :  \la A_1, \cdots , A_n \ra \to \lam_1 A_1 \cup \cdots \cup\lam_n A_n \qv \la \ra : \la \ra\to \emptyset 
\eeq
Note that the underlying discrete opfibration (of categories) $\ul{\widehat\M} \to \C$ 
corresponds to the action of $\C$ on the {\em set} $\P X$
and is the usual symmetry groupoid: $\lam : A \to \lam A$.

\subsection{Tangram game}
\label{ex3}

We can modify the above example by defining the monoid $(\P X , \cup' ,\emptyset)$ 
where $A\cup' B = A\cup B$ if $A$ and $B$ are disjoint, and $A\cup' B = X$ if they are not.
We so obtain the multicategory $\widehat\M'$ with figures as objects and multiarrows
\eq   \label{}
\la \ldots \ra :  \la A_1, \cdots , A_n \ra \to A \qv  \la \ra : \la \ra \to \emptyset 
\eeq
where $A$ is a {\em disjoint} union $\lam_1 A_1 \cup \cdots \cup\lam_n A_n$ 
(if $A$ is not $X$).
Passing to the associated multiposet we get $\la A_1, \cdots , A_n \ra \vdash A$
iff $A$ can be obtained by combining the $A_i$ as in the Tangram game.

\subsection{Coverings}
\label{ex4}

Following example \ref{ex2} above, 
we can exploit the posetal structure of the monoid $(\P X , \supseteq ;\cup ,\emptyset)$,
making it a non-discrete monoidal category. 
Now $\la \ldots \ra :  \la A_1, \cdots , A_n \ra \to A$ is an arrow in $\widehat\M$ 
iff $A \subseteq \lam_1 A_1 \cup \cdots \cup\lam_n A_n$.

Passing to the associated multiposet we get $\la A_1, \cdots , A_n \ra \vdash A$
iff $A$ can be covered by figures congruent with the $A_i$.



\section{Conclusions}

We have seen that, within cartesian multicategories, 
an fp-functor $\C\t \to \D\t$ is just an additive functor, 
an fp-functor $\C_\tm \to \D_\tm$ is just a finite product preserving functor and 
an fp-functor $\C\t \to \Set_\tm$ (more generally $\C\t \to \D_\tm$ or $\C\t \to \M$) 
is just a $\C$-module.

Among algebraic theories (cartesian multicategories),
those of the form $\C\t$ for a preadditive $\C$ thus occupy a special place.
Their models (generalized modules) in any cartesian multicategory $\M$
correspond to additive functors $\C \to \cMon(\M)$.
The (logically) non-linear axiom 
\[
(\lam + \mu)x = \lam x + \mu x 
\]
can be seen either as expressing the fact that $\C\t \to \M$ preserves the cartesian structure
or as expressing the additivity of $\C \to \cMon(\M)$.

Similarly, on the plain level, among \"linear" or \"strongly regular" theories (multicategories)
those of the form $\C\t$ for a category $\C$ occupy a special place.
Their models in any multicategory $\M$ correspond to functors $\C \to \Mon(\M)$.

Finally, we have shown how natural instances of
(promonoidal) multicategories arise as (the domain of)
the fibration corresponding to a (strong) indexed monoidal category.

\newpage



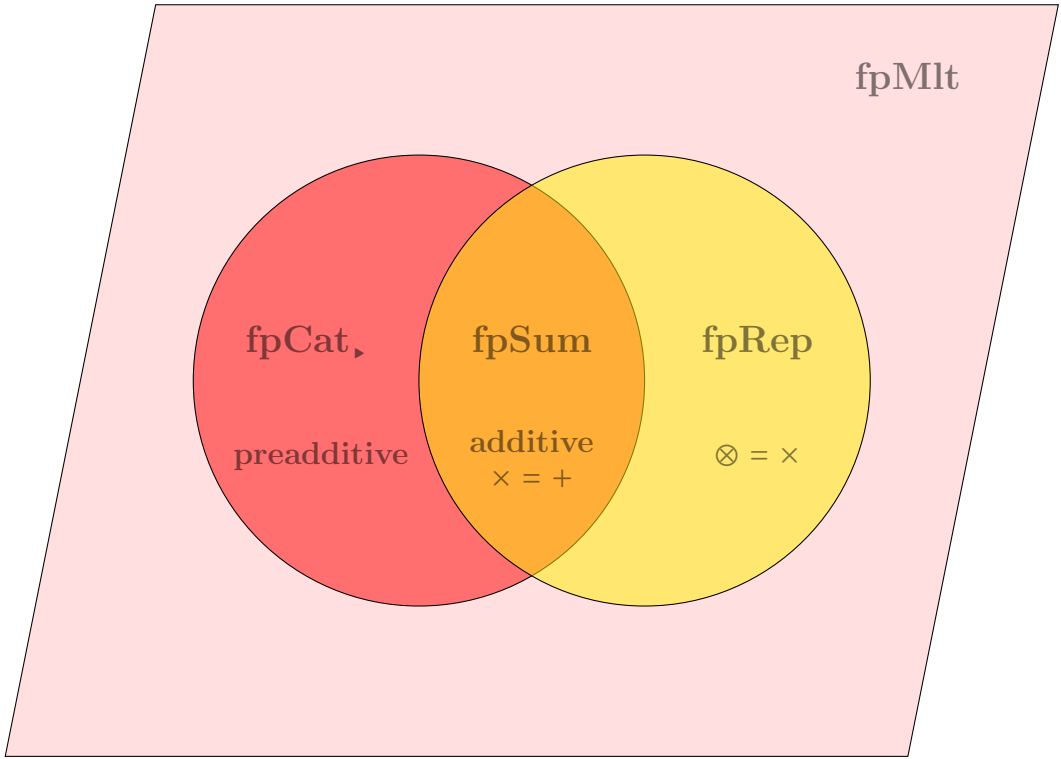
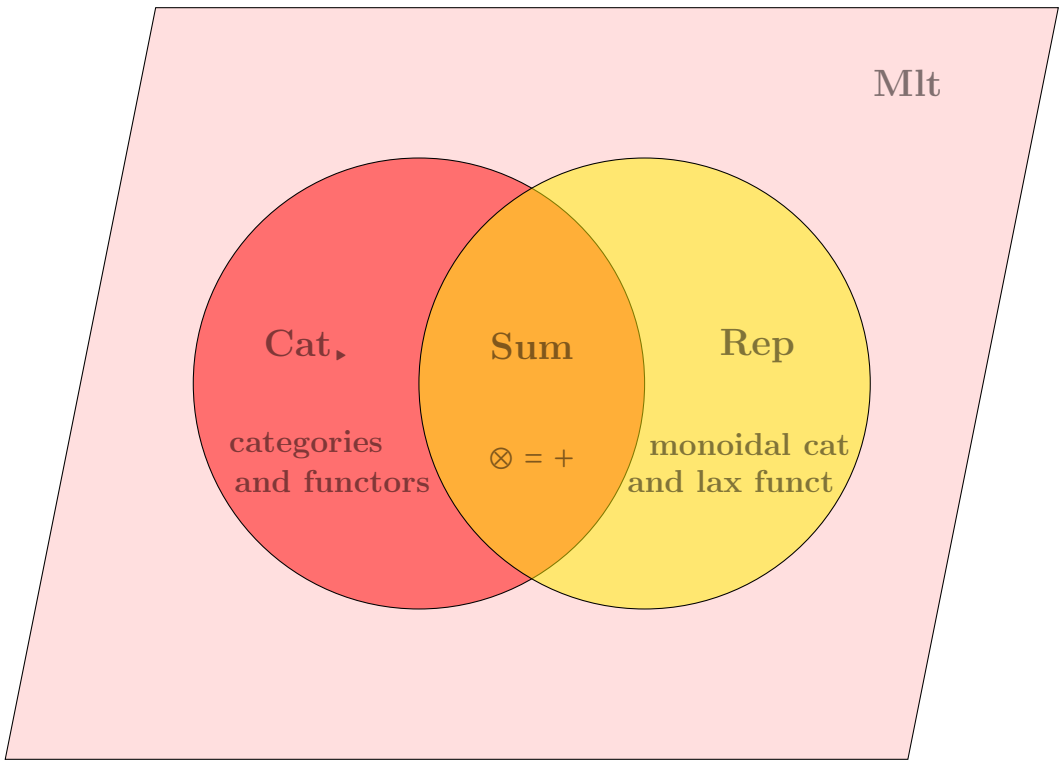
\begin{tikzpicture}
    \begin{scope}[shift={(3cm,-5cm)}, fill opacity=0.5]
        \draw[fill=pink, draw = black] (-5,6) -- (7,6) -- (5,-4) -- (-7,-4) -- cycle;
        \draw[fill=red, draw = black] (-1.5,1) circle (3);
    \draw[fill=yellow, draw = black] (1.5,1) circle (3);
    \node at (5,5) {\large\textbf{\fpMlt}};
    \node at (-3,1.5) {\large\textbf{$\fpCat\t$}};
    \node at (-2.8,0) {\textbf{preadditive}};
    \node at (3,1.5) {\large\textbf{\fpRep}};
    \node at (3,0) {\large\textbf{$\otm = \tm$}};
    \node at (0,1.5) {\large\textbf{\fpSum}};
    \node at (0,0.2) {\textbf{additive}};
    \node at (0,-0.3) {\large\textbf{$\tm = +$}};
    \end{scope}

    \begin{scope}[shift={(3cm,-17cm)}, fill opacity=0.5]
        \draw[fill=pink, draw = black] (-5,6) -- (7,6) -- (5,-4) -- (-7,-4) -- cycle;
        \draw[fill=red, draw = black] (-1.5,1) circle (3);
    \draw[fill=yellow, draw = black] (1.5,1) circle (3);
    \node at (5,5) {\large\textbf{\Mlt}};
    \node at (-3,1.5) {\large\textbf{$\Cat\t$}};
    \node at (-3,0.2) {\textbf{categories}};
    \node at (-2.65,-0.3) {\textbf{and functors}};
    \node at (3,1.5) {\large\textbf{\Rep}};
    \node at (2.9,0.2) {\textbf{monoidal cat}};
    \node at (2.65,-0.3) {\textbf{and lax funct}};
    \node at (0,1.5) {\large\textbf{\Sum}};
    \node at (0,0) {\large\textbf{$\otm = +$}};
    \end{scope}

\end{tikzpicture}


\begin{refs}

\bibitem[Day et al., 2005]{day} B. Day, E. Panchadcharam, R. Street (2005), {\em On centres and lax centres for promonoidal categories},
Colloque International  Charles Ehresmann : 100 ans Universite de Picardie Jules Verne, Amiens.

\bibitem[M. Fiore, 2007]{fiore} M. Fiore (2007) {\em Differential structure in models of multiplicative biadditive intuitionistic linear logic} 
in: Typed Lambda Calculi and Applications, Lecture Notes in Computer Science Volume 4583, 163-177.

\bibitem[T. Fiore, 2005]{fiore2} T. Fiore (2006), {\em Pseudo limits, biadjoints, and pseudo algebras : 
Categorical foundations of conformal field theory}, Memoirs of the American Mathematical Society, 182(860), math.CT/04028298.

\bibitem[Fox, 1976]{fox} T. Fox (1976) Coalgebras and cartesian categories, {\em Communications in Algebra}, {\bf 4(7)}, 665-667.

\bibitem[Gould, 2008]{gould} M. Gould (2008), {\em Coherence for categorified operadic theories}, PhD Thesis, math.CT/1002.0879.

\bibitem[Hermida, 2000]{hermida} C. Hermida (2000), Representable multicategories, {\em Advances in Math.}, {\bf 151}, 164-225.

\bibitem[Hermida, 2004]{hermida2} C. Hermida (2004), {\em Fibrations for abstract multicategories}, in: 
Galois Theory, Hopf Algebras and Semiabelian Categories, Fields inst. comm. AMS, 281-293.

\bibitem[Heunen, 2008]{heunen} C. Heunen (2008), Semimodule enrichment, {\em Electr. Notes in Th. Comp. Sci.}, {\bf 218}, 193-208.

\bibitem[Lack, 2009]{lack} S. Lack (2009), {\em Non-canonical isomorphisms}, math.CT/0912.2126.

\bibitem[Lamarche, 2007]{lamarche} F. Lamarche (2007), Exploring the gap between linear and classical logic,
{\em Theory and Appl. Cat.} {\bf 18}, 471-533.

\bibitem[Lambek, 1989]{lambek} J. Lambek (1989), {\em Multicategories revisited}, 
in: Contemporary Mathematics 92 (Amer. Math. Soc., Providence), 217-239.


\bibitem[Leinster, 2003]{leinster} T. Leinster (2003), {\em Higher operads, higher categories}, Cambridge University Press, math.CT/0305049.

\bibitem[Pavlovic, 2008]{pavlovic} D. Pavlovic (2008), {\em Geometry of abstraction in quantum computation}, in:
Proceedings of Symposia in Applied Mathematics, AMS, {\bf 71}, 233-267.

\bibitem[Shulman, 2008]{shulman} M. Shulman (2008), Framed bicategories and monoidal fibrations, 
{\em Theory and Appl. Cat.} {\bf 20}, 650-738.

\end{refs}

\end{document}